\documentclass[12pt]{article}
\usepackage{fullpage}
\usepackage{graphicx}
\usepackage{amsmath}
\usepackage{amssymb}

\newtheorem{theorem}{Theorem}
\newtheorem{lemma}[theorem]{Lemma}
\newtheorem{corollary}[theorem]{Corollary}

\newcommand{\R}{\mathbb{R}}
\newcommand{\T}{\mathbb{T}}
\newcommand{\Z}{\mathbb{Z}}

\newcommand{\cP}{\mathcal{P}} 

\newcommand{\proof}{\noindent {\bf Proof: }}
\newcommand{\qed}{\hfill {\bf QED}}

\begin{document}
\title{Multi-cluster dynamics in coupled phase oscillator networks}
\author{Asma Ismail\thanks{Also: Department of Mathematics, University of Benghazi, Benghazi, Libya} and Peter Ashwin\\
Mathematics Research Institute,\\
Harrison Building,\\
University of Exeter,\\Exeter EX4 4QF, UK
}

\maketitle

\begin{abstract}
In this paper we examine robust clustering behaviour with multiple nontrivial clusters for identically and globally coupled phase oscillators. These systems are such that the dynamics is completely determined by the number of oscillators $N$ and a single scalar function $g(\varphi)$ (the coupling function). Previous work has shown that (a) any clustering can stably appear via choice of a suitable coupling function and (b) open sets of coupling functions can generate heteroclinic network attractors between cluster states of saddle type, though there seem to be no examples where saddles with more than two nontrivial clusters are involved. In this work we clarify the relationship between the coupling function and the dynamics. We focus on cases where the clusters are inequivalent in the sense of not being related by a temporal symmetry, and demonstrate that there are coupling functions that give robust heteroclinic networks between periodic states involving three or more nontrivial clusters. We consider an example for $N=6$ oscillators where the clustering is into three inequivalent clusters. We also discuss some aspects of the bifurcation structure for periodic multi-cluster states and show that the transverse stability of inequivalent clusters can, to a large extent, be varied independently of the tangential stability.
\end{abstract}


\section{Introduction}

Coupled oscillator models exhibit complex dynamics that has been observed in a wide range of different fields including physical \cite{Aronson-Krupa-Ashwin96,Gonzales-Masoller-Torrent-Garcia2007} and biological \cite{Mikhail-Huerta-Varona-Afraimovich2006} models. Synchronization \cite{Strogatz2000}, clustering \cite{Okuda93,Brown-Holmes-Moehlis2003,Kori-etal2014}, chaos \cite{Bick-et-al} and spontaneous switching between different cluster states \cite{Orosz-Ashwin-Townley2009,Kori-etal2014} have all been observed in such systems. Other studies have examined coupling between two or more systems that may individually be chaotic, and a wide variety of types of synchronization have been found and analysed; see for example \cite{Fujisaka_Yamada83, Pecora_Carrol1990,Vladimi_Igor_Hasler2004}.

We examine phase oscillator models that are appropriate if the coupling between oscillators is weak compared to the attraction onto the limit cycle (e.g. \cite{Ashwin-Swift1992,Brown-Holmes-Moehlis2003}). Although the coupling structure and strength are important for the dynamical behaviour of the system, the exact coupling function (which represents the nonlinearities in the oscillators and the coupling) has a subtle effect on the collective behaviour of the system.

Research into the dynamics of coupled nonlinear oscillators has long explored the question ``What is the dynamics of a given system?''. A less frequently asked, but also very interesting, question is ``How can we design a coupled system to have specific dynamics?''. This latter question was considered by \cite{Orosz-Moehlis-Ashwin2009} who designed cluster states with a prescribed clustering by giving explicit conditions on the coupling function and its first derivative to have a stable cluster state with a specific clustering. They demonstrate specific coupling functions that give stable cluster states for any partition of the oscillator into groups, regardless the number of oscillators and the size of each cluster.

In this paper we go beyond \cite{Orosz-Moehlis-Ashwin2009} in four ways. Firstly, we examine three-cluster states and show that not only stable cluster states, but also cluster states with specific transverse stability properties can be designed by suitable choice of coupling function.  Secondly, we give some results on how the transverse stability can be varied independently of the tangential stability and hence exhibit possible bifurcation scenarios from transversely stable clustering. Thirdly, we show examples of how nontrivial cluster states with three inequivalent clusters can be joined into a heteroclinic network. Finally, we generalize some of the bifurcation results to more general multi-cluster states with an arbitrary numbers of clusters.

We use a Fourier representation of the coupling function associated with a system of $N$ oscillators to design general three-cluster states, as in \cite{Orosz-Moehlis-Ashwin2009}. The rest of the paper is organized as follows; for the remainder of this section we recall some of the notation and previous results on existence and stability of periodic cluster states. We define a notion of inequivalence of clusters within a cluster state and consider some sufficient conditions for clusters to be inequivalent. Section~\ref{sec:3cluster} recalls and extends some basic results on the appearance of tangentially stable but transversely unstable three-cluster states. We present in Theorem~\ref{thm:transstab} a characterization of transverse stability, and in Corollary~\ref{cor:couplingfunction} a result on transverse bifurcation of three-cluster states. Section~\ref{sec:heteroclinic} presents what we claim is the smallest possible cluster state with three inequivalent nontrivial clusters (requiring at least $N=6$ oscillators) and gives some examples how these may be connected into robust attracting heteroclinic networks \cite{Krupa97}. Finally, Section~\ref{sec:conclusion} discusses some consequences of this work, including a generalization of Corollary~\ref{cor:couplingfunction}.

\subsection{Periodic cluster states and symmetries for globally coupled oscillators}

In this paper we consider $N$ phase oscillators that are all-to-all coupled and governed by the following generalization \cite{Ashwin-Swift1992,Okuda93,Hiroaki96,Orosz-Moehlis-Ashwin2009} of Kuramoto's model system of coupled phase oscillators:
\begin{equation}
\frac{d\theta_i}{dt}=\omega+\frac{1}{N}\sum_{j=1}^N g(\theta_i-\theta_j)
\label{eq:coupledoscillatorsystem}
\end{equation}
where $\theta_i\in \T=[0,2\pi)$ is the phase of the $i^{th}$ oscillator, $i=1,\cdots,N$ and $g:\T\rightarrow \R$ is a $2\pi$-periodic nonlinear {\em coupling function} that we assume is smooth and represented by a truncated Fourier series as in  \cite{Hiroaki96,Okuda93}:
\begin{equation}
g(\varphi)=\sum_{r=0}^R\left(c_r \cos(r\varphi)+s_r \sin(r\varphi)\right)
\label{eq:couplingfunction}
\end{equation}
where $c_r$ ($r\geq 0$) and $s_r$ ($r>0$) are the real coefficients and $R$ is the number of Fourier modes. Note that the coupling function $g(\varphi)$ derived from weakly coupled nonlinear phase oscillators will typically have several non-zero modes in its Fourier series, even if the oscillators are close to Hopf bifurcation \cite{Kori-etal2014}.

Conditions on the coupling function and its first derivative that ensure the existence and stability of desired cluster states in the system (\ref{eq:coupledoscillatorsystem}) are derived in \cite{Orosz-Moehlis-Ashwin2009}. Note that the system is invariant under ``spatial'' symmetries $S_N$ acting by permutation of the components and ``temporal'' symmetries $\T$ given by 
\begin{equation}
(\theta_1,\ldots,\theta_N)\mapsto (\theta_1+\phi,\ldots,\theta_N+\phi)
\label{eq:Tsymm}
\end{equation}
for any $\phi\in\T$ \cite{Ashwin-Swift1992}. We now look at periodic cluster states in a bit more depth. Consider a partition $\cP=\{p_1,p_2,...,p_M \}$ into $M$ clusters, where $1 \leq M \leq N$; each of $p_k$ form a cluster of size $m_k=\left|p_k\right|$ for $k=1,...,M$ and $\sum_{k=1}^M m_k =N$.  We say a cluster is a {\em multi-cluster} if $M\geq 3$. The $k$th cluster is said to be {\em nontrivial} if $m_k>1$.

A periodic orbit $\theta(t)=(\theta_1(t),\ldots,\theta_N(t))$ of (\ref{eq:coupledoscillatorsystem}) defines an associated {\em clustering} $\cP[\theta]$ where $i$ and $j$ being the same cluster $p_k$ if and only if $\theta_i(t)=\theta_j(t)$ for some $t$. It is possible that periodic cluster states can have additional symmetries associated with non-trivial phase shift symmetry $\T$; more precisely, it is known that they can be classified according to decompositions of the form
$$
N=m(k_1+\cdots+k_{\ell})
$$
where $m,k_i,\ell$ are all positive integers.

\begin{theorem}[\cite{Ashwin-Swift1992}, Thm~3.1]\label{thm:ASiso}
The subsets of $\T^N$ that are invariant for (\ref{eq:coupledoscillatorsystem}) because of symmetries of $S_N\times \T$ corresponding to isotropy subgroups in the conjugacy class of
$$
\Sigma_{k,m}:=(S_{k_1}\times\cdots\times S_{k_{\ell}})^m \times_s \Z_m
$$
where $N=mk$, $k=k_1+\cdots+k_{\ell}$ and $\times_s$ denotes the semidirect product. The points with this isotropy have $\ell m$ clusters that are split into $\ell$ groups of $m$ clusters of the size $k_i$. The clusters within these groups are cyclically permuted by a phase shift of $2\pi/m$. The number of isotropy subgroups in this conjugacy class is $N!/[m(k_1!\ldots k_\ell!)]$.
\end{theorem}

We say two clusters of $\cP[\theta]$ are {\em equivalent} if there is a symmetry in $S_N\times \T$ that maps one cluster to the other. Otherwise they are said to be {\em inequivalent}. A clustering is said to be {\em phase non-degenerate} if the phase difference between two different clusters is only attained by those clusters; more precisely we say the cluster phases $\{\phi_k\}_{k=1}^M$ are {\em non-degenerate} if 
\begin{equation}
d(\phi_i-\phi_j,\phi_k-\phi_\ell)>0
\label{eq:phasenondeg}
\end{equation}
for all $i,j,k,\ell\in\{1,\ldots,M\}$ with $i\neq $j, $k\neq \ell$, unless $i=k$ and $j=\ell$.\footnote{We choose $d(\phi,\psi)=1-\cos(\phi-\psi)$ as a metric on the circle.} As noted in \cite{Orosz-Moehlis-Ashwin2009}, a sufficient condition for this is that the phase differences are rationally independent of each other and of $2\pi$. The following more general statement follows from the definitions and Theorem~\ref{thm:ASiso}; note that if the phase differences are rationally independent this implies (c).

\begin{lemma}
Suppose $\theta\in \T^N$ is a cluster state. Any of the following is a sufficient condition for all clusters to be inequivalent:
\begin{itemize}
\item[(a)] $N$ is prime and at least one cluster is nontrivial.
\item[(b)] no two clusters have the same size
\item[(c)] the clustering is phase non-degenerate
\end{itemize}
\end{lemma}

\proof
In case (a), Theorem~\ref{thm:ASiso} gives that the only possible factorizations $N=m(k_1+\cdots+k_{\ell})$ have $m=1$ or $m=N$. The latter case is ruled out if one of the $k_i>1$. In case (b), this also implies that $m=1$ in Theorem~\ref{thm:ASiso}. Finally in case (c), if $m>1$ then there will be clusters which can be chosen such that $\phi_1-\phi_2=2\pi/m$ and $\phi_2-\phi_k=2\pi/m$ ($k=1$ if $m=2$, otherwise we can choose $k=3$. This implies that the clustering will be phase degenerate and conversely, phase non-degeneracy implies that $m=1$ and the clusters are inequivalent.
\qed

~~

Given a partition $\cP$ let us now define the subspace $\T_\cP^N\subset \T^N$ in which clustering occurs as follows:
\begin{center}
$\T_\cP^N$=$\{ \theta \in \T^N$: If there exists $k$ such that $i,j\subset p_k$ then $\theta_i=\theta_j\}.$
\end{center}
Because (\ref{eq:coupledoscillatorsystem}) is equivariant under the action of $S_N$ by permutation of the oscillators and by the circle group $\T$ (\ref{eq:Tsymm}), these cluster states $\T_\cP^N$ are simply fixed point subspaces for groups conjugate to $S_{p_1}\times S_{p_2}\times \cdots \times S_{p_M}$, i.e. with $m=1$.

The system (\ref{eq:coupledoscillatorsystem}) can be simplified on $\T_\cP^N$ as follows; as if $i$ is in the 
$k$th cluster then we say $\theta_i=\Psi_k$ and the system reduces to
\begin{equation}
\frac{d\Psi_k}{dt}= \omega + \frac{1}{N} \sum_{\ell=1}^M m_{\ell} g(\Psi_k-\Psi_{\ell}),
\label{eq:phasemodel}
\end{equation}
for $k=1,\ldots,M$.

The dynamics of a periodic multi-cluster state can then be expressed as $\Psi_k= \phi_k+\Omega t$ for $k=1,...,M$ where $\phi_k\in \T$ represents the relative phases of clustering and $\Omega \in \R^+$ is the frequency of the periodic orbit. This means that:
\begin{equation}
\Omega= \omega+ \frac{1}{N} \sum_{\ell=1}^M m_{\ell} g(\phi_k-\phi_{\ell})
\label{eq:frequency}
\end{equation}
where $k=1,...,M$. Computing the difference between the first equation and the remaining ones gives:
$\sum_{\ell=1}^M m_{\ell}\left(g(\phi_k- \phi_{\ell})- g(\phi_1- \phi_{\ell})\right)=0$. Defining 
$$
g_0:=g(\phi_k- \phi_k)\mbox{ and } g_{k\ell}:=g(\phi_k- \phi_{\ell})~~~(k\neq \ell)
$$
as in \cite{Orosz-Moehlis-Ashwin2009} means that we can rewrite (\ref{eq:frequency}) as
\begin{equation}
\sum_{\ell=1}^M m_{\ell}(g_{k\ell}-g_{1\ell})=0,~~k=2,...,M.
\label{eq:existenceequations}
\end{equation}
Equations (\ref{eq:existenceequations}) are conditions that that the coupling function should satisfy for existence of such a periodic cluster state. We also refer to \cite{Okuda93,Orosz-Moehlis-Ashwin2009} for a discussion of stability where it is shown that a cluster state is linearly stable if and only if it is linearly stable in both of the following senses:

\begin{itemize}
\item {\bf Tangential stability:} (also called inter-cluster stability \cite{Okuda93}) To determine the stability of the state to the change in its phases that respects the clustering we consider the linearized stability for perturbations $\Psi_k =\omega t+\phi_k$, and writing 
$\chi=\Psi-(\phi_k +\Omega t)$, gives the following:
\begin{equation}
\frac{d\chi_k}{dt} = \frac{1}{N} \sum_{\ell=1}^M m_{\ell} g'(\phi_k - \phi_{\ell})(\chi_k-\chi_{\ell})=\sum_{{\ell}=1}^M T_{k,{\ell}}\chi_{\ell}
\end{equation}
where, as stated in \cite{Orosz-Moehlis-Ashwin2009}, $T$ is the matrix:
\[
T_{k,{\ell}}=\frac{1}{N}\left[\delta_{k,{\ell}}\left(\sum_{r=1,r\neq k}^M m_r g'_{k,r}\right)-(1-\delta_{k,{\ell}}m_{\ell} g'_{k,{\ell}})\right]
\]
where $\delta$ is the Kronecker delta. The matrix $T$ has $M$ eigenvalues (including one trivial value) and for tangential stability we require that all the other tangent eigenvalues have negative real parts. This means that a $p^{th}$ cluster is tangentially stable if:
\begin{equation}
Re(\lambda_k^{tang})<0,~~~~k=2,...,M.
\label{eq:tangentstability}
\end{equation}
\item {\bf Transverse stability:} (also called intra-cluster stability \cite{Okuda93}) The stability of a periodic state to changes of phases that change the clustering and is obtained by linearizing (\ref{eq:coupledoscillatorsystem}) about $\eta_i=\theta_i-(\phi_k+\Omega t)$, which gives the following :
\[
\frac{d\eta_u}{dt}=\sum_{v=1}^{s_k} S_{u,v}^{(k)}~\eta_v, 
\]
where (as stated in \cite{Orosz-Moehlis-Ashwin2009}) $S$ is given by 
\[ 
S_{u,v}^{(k)}= \frac{1}{N}\left[\delta_{u,v}\left(\sum_{r=1}^M m_r g'_{k,r}-g'_0\right)-\left(1-\delta_{u,v}\right)g'_0\right].
\]
This has $(M+k)$ real eigenvalues that are negative for a transverse stable cluster state; put otherwise, a cluster state is transversely stable  if:
 \begin{equation}
 \lambda_{M+k}^{tran}<0,~~~k=1,...,W.
 \label{eq:transversestability}
 \end{equation}
 where $W$ is the number of nontrivial clusters (i.e. those with more than one oscillator).
\end{itemize} 
A coupling function associated with a stable cluster state will satisfy (\ref{eq:existenceequations}), (\ref{eq:tangentstability}), and (\ref{eq:transversestability}) and we observe that the number of degrees of freedom on choosing the $g(\phi)$ will always be surplus to requirements \cite{Orosz-Moehlis-Ashwin2009} if enough Fourier modes are chosen. Note that the above holds for cluster states regardless of whether the clusters are inequivalent or not. If two clusters are equivalent then the transverse eigenvalues for those clusters will be equal; if the clusters are inequivalent then generically the transverse eigenvalues for those clusters will be unequal.

\section{Three-cluster states and their stability}
\label{sec:3cluster}

We now derive conditions for tangential and transverse stability of three-cluster states for the system  (\ref{eq:coupledoscillatorsystem}) of globally coupled phase oscillators. Consider a periodic state where all clusters are nontrivial (i.e. the clusters have sizes $m_r\geq 2$ for $r=1,2,3$, where $m_1+m_2+m_3=N$). This implies that $N\geq 6$; we explore the special case $N=6$ in more detail in Section~\ref{subsec:N=6}. Recall from \cite{Orosz-Moehlis-Ashwin2009} that the condition for existence of a three-cluster state is:
\begin{equation}
\begin{split}
m_1(g_{21}-g_0)+m_2(g_0-g_{12})+m_3(g_{23}-g_{13})&=0\\
m_1(g_{31}-g_0)+m_2(g_{32}-g_{12})+m_3(g_0-g_{13})&=0.
\end{split}
\label{eq:existenceconditionsof3-cluster}
\end{equation}
The tangential stability is determined by (\ref{eq:tangentstability}), namely
\begin{equation}
\begin{split}
\lambda_1^{tang}=&0\\
\lambda_2^{tang}=&\frac{1}{2}(\mu+i\sqrt{\nu-\mu^2}) \\
\lambda_3^{tang}=&\frac{1}{2}(\mu-i\sqrt{\nu-\mu^2})
\end{split}
\label{eq:tangenteigenvalues2004}
\end{equation}
where
\begin{equation}
\mu=\frac{1}{N}(m_2 g'_{12}+m_3 g'_{13}+ m_1 g'_{21}+ m_3g'_{23}+ m_1g'_{31}+ m_2g'_{32}).
\label{eq:mueqn}
\end{equation}
and
\begin{equation}
\begin{split}
\nu=&\frac{4}{N^2}\left((m_1g'_{21} +m_3 g'_{23})(m_1g'_{31} + m_2 g'_{32} \right)- m_2 m_3g'_{32}g'_{23}\\
& +(m_2g'_{12}+m_3g'_{13})(m_1g'_{31}+m_2g'_{32}) -m_1 m_3 g'_{31}g'_{13}\\
& +(m_2g'_{12}+ m_3 g'_{13})(m_1g'_{21}+ m_3 g'_{23}) - m_1 m_2 g'_{21}g'_{12}).
\end{split}
\label{eq:oldnu}
\end{equation}
Transverse stability is determined by the eigenvalues:
\begin{equation}
\begin{split}
\lambda_{4}^{tran}=&\frac{1}{N} (m_1g'_{0}+m_2 g'_{12}+m_3g'_{13})\\
\lambda_{5}^{tran}=&\frac{1}{N} (m_1g'_{21}+m_2 g'_{0}+m_3g'_{23})\\
\lambda_{6}^{tran}=&\frac{1}{N} (m_1g'_{31}+m_2 g'_{32}+m_3g'_{0})
\end{split}
\label{eq:transeigs3}
\end{equation}
where the multiplicities of $\lambda_{4,5,6}^{tran}$ are $m_{1,2,3}-1$ respectively. Our first new result is the following sufficient condition for tangential stability of three-clusters:

\begin{lemma}
Suppose there is a periodic three-cluster state such that $g'_{ij}<0$  for all $i\neq j$. Then the cluster is tangentially stable with complex contracting eigenvalues.
\end{lemma}

\proof
If $g'_{ij}<0$ then all terms in (\ref{eq:mueqn}) are negative and so $\mu<0$. Moreover, note that (\ref{eq:oldnu}) can be written in the form
\begin{equation}
\begin{split}
\nu= &\frac{4}{N^2} \left[m_1^2g'_{21}g'_{31} +m_2^2g'_{12}g'_{32}+m_3^2 g'_{13}g'_{23}\right.\\
 &+ m_1 m_2\left(g'_{21}g'_{32}+ g'_{12}g'_{31}\right) +m_1m_3 \left(g'_{23}g'_{31}+ g'_{13}g'_{21}\right) \left.+m_2m_3\left(g'_{13}g'_{32}+g'_{12}g'_{23}\right) \right].
\end{split}
\label{eq:newnu}
\end{equation}
Hence if $g'_{ij}<0$ then all terms in (\ref{eq:newnu}) are positive and so $\nu>0$. Hence the eigenvalues (\ref{eq:tangenteigenvalues2004}) are complex with negative real parts and so the cluster is tangentially stable.
\qed

~

The next result demonstrates that the tangential and transverse stability can be set independently of each other. We define
$$
K_1:=\frac{1}{m_1} (m_2g'_{12}+m_3g'_{13}),~~
K_2:=\frac{1}{m_2} (m_1g'_{21}+m_3g'_{23}),~~
K_3:=\frac{1}{m_3} (m_1g'_{31}+m_2g'_{32})
$$
Without loss of generality (renumbering the clusters if necessary) let us assume that 
$$
K_1\leq K_2\leq K_3
$$
and we demonstrate the following:

\begin{theorem}
\label{thm:transstab}
Suppose that $g(\varphi)$ is such that there is a periodic three-cluster state with non-trivial clusters such that the clusters are tangentially stable and assume that $K_1\leq K_2\leq K_3$. Then we can classify the transverse stability as follows:
\begin{itemize}
\item If $-g'_0<K_1$ then $\lambda_{4,5,6}^{tran}>0$ (all clusters unstable).
\item If $K_1<-g'_0<K_2$ then $\lambda_{4}^{tran}<0$ and $\lambda_{5,6}^{tran}>0$ (one stable cluster).
\item If $K_2<-g'_0<K_3$ then $\lambda_{4,5}^{tran}<0$ and $\lambda_{6}^{tran}>0$  (two stable clusters).
\item If $K_3<-g'_0$ then $\lambda_{4,5,6}^{tran}<0$ (all clusters stable).
\end{itemize}
\end{theorem}

\proof
These conclusions follows from noting that the conditions on $g'_0$ ensure that (\ref{eq:transeigs3}) have zero, one, two or three positive transverse eigenvalues as stated. Note that this is independent of the number of oscillators, although if one or more of the clusters are trivial, the transverse exponent for that cluster is not defined. If $-g'_0=K_i$ for any $i$ then the state will be at a bifurcation point and the stability is not determined at linear order.
\qed

~

Theorem~\ref{thm:transstab} can be used to prove the following Corollary about bifurcation of three-cluster states involving changes in transverse stability.

\begin{corollary}
\label{cor:couplingfunction}
Suppose that $g(\varphi)$ is such that there is a periodic three-cluster state with non-trivial clusters, such that the clusters are tangentially stable. Then there is a parametrized family of coupling functions
$$
g_{r}(\varphi)=g(\varphi)+r h(\varphi)
$$
and parameter values $r_1 \leq r_2 \leq r_3$ such that for all values of the parameter $r\in\R$ the cluster state remains with the same phases and tangentially stability and:
\begin{enumerate}
\item If $r<r_1$ then $\lambda_{4,5,6}^{tran}>0$ (all clusters unstable).
\item If $r_1<r<r_2$ then $\lambda_{4}^{tran}<0$ and $\lambda_{5,6}^{tran}>0$ (one stable cluster).
\item If $r_2<r<r_3$ then $\lambda_{4,5}^{tran}<0$ and $\lambda_{6}^{tran}>0$  (two stable clusters).
\item If $r_3<r$ then $\lambda_{4,5,6}^{tran}<0$ (all clusters stable).
\end{enumerate}
\end{corollary}

\proof
For this specific three cluster state with relative phases $(\phi_1,\phi_2,\phi_3)$ there will be a $0<\epsilon<\pi$ such that $d(\phi_j-\phi_k,0)>\epsilon$ for all $j\neq k$. Now consider a smooth compactly supported periodic function $h$ such that $h(\varphi)=0$ for all $\varphi$ with $d(\varphi,0)>\epsilon$, $h(0)=0$ and $h'(0)=-1$. One can verify that $(g_{r})_{i,j}=g_{i,j}$, $(g_{r})'_{i,j}=g'_{i,j}$ for all $i\neq j$ and $r\in\R$, $g_{r}(0)=g(0)$ and 
$$
g'_{r;0}=g'_0-r.
$$
for all $r\in \R$. Hence the existence condition (\ref{eq:existenceequations}) and the tangential stability conditions do not depend on $r$, while the cases of Theorem~\ref{thm:transstab} translate into the cases depending on $r$.
\qed

~

\subsection{Cluster dynamics for $N=6$ oscillators}
\label{subsec:N=6}

In this section we consider properties of cluster states with nontrivial and inequivalent clusters for $N=6$, before applying the results from the previous section to give sufficient conditions for a nontrivial stable three cluster state for (\ref{eq:coupledoscillatorsystem}). By restricting to phase non-degenerate (and hence inequivalent) clusters we do not consider the cases $m>1$ in Theorem~\ref{thm:ASiso}. It does not appear to be easy to characterise the solutions of the system of equation for the Fourier coefficients analytically. For $N=6$ recall that (\ref{eq:coupledoscillatorsystem}) can be written as
\begin{equation}
\frac{d\theta_i}{dt}=\omega+\frac{1}{6}\sum_{j=1}^6 g(\theta_i-\theta_j).
\label{eq:sixcoupledoscillatorsystem}
\end{equation}
The symmetries of interchange of the $6$ oscillators, $S_6$, gives nine possible isotropy subgroups corresponding to possible cluster states that are listed in the Table \ref{tab:FixedpointsubspacesofS_6}; each periodic cluster state will reside in precisely one of these invariant subspaces.

\begin{table}
\begin{center} 
\begin{tabular}{|c|c|c|r|r|}
\hline
Isotropy & dim(Fix($\Sigma$))& Representative & Number of& Orbit size\\
subgroup $\Sigma$ & & point &  conjugates & \\
\hline
$S_6$&1&$(\theta_1,\theta_1,\theta_1,\theta_1,\theta_1,\theta_1)$&1&1\\
\hline
$S_5$&2&$(\theta_1,\theta_1,\theta_1,\theta_1,\theta_1,\theta_2)$&6&6\\
\hline
$(S_3)^2$&2 &$(\theta_1,\theta_1,\theta_1,\theta_2,\theta_2,\theta_2)$&10&20\\
\hline
$S_4\times S_2$&2&$(\theta_1,\theta_1,\theta_1,\theta_1,\theta_2,\theta_2)$&15&15\\
\hline
$S_4$ &3&$(\theta_1,\theta_1,\theta_1,\theta_1,\theta_2,\theta_3)$&30&60\\
\hline
$S_2\times S_3$&3&$(\theta_1,\theta_1,\theta_2,\theta_2,\theta_2,\theta_3)$&60&60\\
\hline
$(S_2)^3$&3&$(\theta_1,\theta_1,\theta_2,\theta_2,\theta_3,\theta_3)$&15&90\\
\hline
$(S_2)^2$&4&$(\theta_1,\theta_1,\theta_2,\theta_2,\theta_3,\theta_4)$&45&180\\
\hline
$S_3$&4&$(\theta_1,\theta_1,\theta_1,\theta_2,\theta_3,\theta_4)$&20&120\\
\hline
$S_2$&5&$(\theta_1,\theta_1,\theta_2,\theta_3,\theta_4,\theta_5)$&15&360\\
\hline
$I$&6&$(\theta_1,\theta_2,\theta_3,\theta_4,\theta_5,\theta_6)$&1&720\\
\hline
\end{tabular}
\end{center}
\caption{Conjugacy classes of isotropy subgroups and representative fixed-points subspaces for the action of $S_6$ on the phase space for $N=6$ globally coupled oscillators corresponding to inequivalent clusters. Note that $I$ represents the trivial group while the number of conjugate groups and the number of point in the group orbit under $S_6$ are given in the last two columns.}
\label{tab:FixedpointsubspacesofS_6} 
\end{table}


To summarise the calculations from the previous section, there is a stable cluster state of this type if we can solve the two equations and five inequalities:
\begin{equation}  
\begin{split}
2(g_{21}-g_0)+2(g_0-g_{12})+2(g_{23}-g_{13})&=0;\\
2(g_{31}-g_0)+2(g_{32}-g_{12})+2(g_0-g_{13})&=0; \\
Re(\lambda_2^{tang})=Re(\frac{1}{2}(\mu+i\sqrt{\nu-\mu^2}))&<0; \\
Re(\lambda_3^{tang})=Re(\frac{1}{2}(\mu-i\sqrt{\nu-\mu^2}))&<0;\\
Re(\lambda_4^{tran})=\frac{1}{3}(g'_0+g'_{12}+g'_{13})&<0;\\
Re(\lambda_5^{tran})=\frac{1}{3}(g'_{21}+g'_0+g'_{23})&<0; \\
Re(\lambda_6^{tran})=\frac{1}{3}(g'_{31}+g'_{32}+g'_0)&<0.
\end{split}
\label{eq:stable3-clusterstateconditions}
\end{equation}
If a cluster is phase non-degenerate then the six values $g_{i,j}$ for $i\neq j$, $g_0$ (and the derivatives at these points) can be chosen independent of each other. This gives fourteen degrees of freedom to satisfy only seven constraints.

Indeed, it is known that (\ref{eq:stable3-clusterstateconditions}) can be satisfied \cite{Orosz-Moehlis-Ashwin2009} by choice of a single high order trigonometric coupling function of the form 
$$
g(\varphi)=-\sin(4\varphi)
$$
for the state $(\psi_1,\psi_2,\psi_3)=(\psi,\psi+\pi/2,\psi+\pi)$ so that $g_{ij}=g_0=0$ and $g'_{ij}=g'_0=-1$; for this $\mu<0$ and $\nu>0$. However, this state is phase degenerate as, for example $g_{12}=g_{23}$. The simplest coupling function that gives a phase non-degenerate stable three-cluster state will be of higher order; for example
$$
g(\varphi)=-\sin(L\varphi)
$$
for $L\in \Z$, $L\geq 5$ will stabilize the phase non-degenerate state $(\psi_1,\psi_2,\psi_3)=(\psi,\psi+2\pi/L,\psi+6\pi/L)$. Applying Corollary~\ref{cor:couplingfunction} to this state shows that there are perturbations of this coupling function that can have a range of transverse stabilities. The next section considers some explicit examples of phase non-degenerate states, and heteroclinic cycles that connect these. 

\section{Heteroclinic attractors involving three-cluster states}
\label{sec:heteroclinic}

Considering $L=4$ and the Fourier coefficients $(c_r,s_r)$ for $r=1,\ldots,L$ we specify the coupling function for the system (\ref{eq:sixcoupledoscillatorsystem}). Numerical investigations of the dynamics reveal the existence of a range of complex dynamics, including attraction to heteroclinic cycles between three-cluster states that have one or two transversely unstable directions. In the following we concentrate on analyzing the dynamics of the heteroclinic cycle that are formed by connections between non-trivial $3-$cluster states. Table~\ref{tab:3clustereg} lists values of these Fourier coefficients.

\begin{table}[ht]
\begin{center}
	\begin{tabular}{|p{2cm}|l|l|l|l|l|l|l|l|}
	\hline
& $c_1$ &$c_2$ &$c_3$ &$c_4$ &$s_1$ &$s_2$ &$s_3$ & $s_4$ \\
\hline
  \textbf{Case 0}&  $0$ & $0$ & $0$ & $0$ & $0$ & $0$ & $0$ & $-1$ \\
	\hline
	\textbf{Case 1} & $0.31185$ & $0.37096$& $0$& $0.99008$ & $0.10793$& $0.58180$& $0$& $-0.14053$\\
	\hline 
	\textbf{Case 2}&  $0.31185$ & $0.39$& $0$& $0.99008$ & $0.10793$& $0.58180$& $0$& $-0.14053$\\
	\hline
	\end{tabular}

\vspace{3mm}

	\begin{tabular}{|l|l|l|l|l|}
	\hline
 & $\alpha=\Psi_1-\Psi_3$ &$\beta=\Psi_2-\Psi_3$ & $\lambda^{tang}_i$ & $\lambda^{tran}_i$ \\
\hline
	 \textbf{Case 0} & $\pi/2$ & $\pi$ & $0$, $-1\pm i$  &  $-1$, $-1$, $-1$ \\
	\hline
	 \textbf{Case 1} & $1.7014$ & $4.7573$ & $0$, $-0.4473$, $-1.4690$  &  $-1.3070$, $-0.06014$, $0.1636$ \\
	\hline 
	\textbf{Case 2} & $1.7087$ & $4.7761$ & $0$, $-0.5102$, $-1.3901$  &  $-1.2798$, $0.03692$, $0.02568$ \\
	\hline
	\end{tabular}
\end{center}	
\caption{\label{tab:3clustereg}
Top: Fourier coefficients giving rise to attractors that include three-cluster states of type $(2,2,2)$. Bottom: properties of the corresponding three-cluster states. Case 0 gives a stable $(2,2,2)$ state while for Cases 1,2 these states are transversely unstable. Note that Figure~\ref{fig:examplecase12} shows that the $(2,2,2)$ state can appear within an attracting heteroclinic network for the given parameter  values.}
\end{table}

\begin{figure}%
\centerline{\includegraphics[width=8cm]{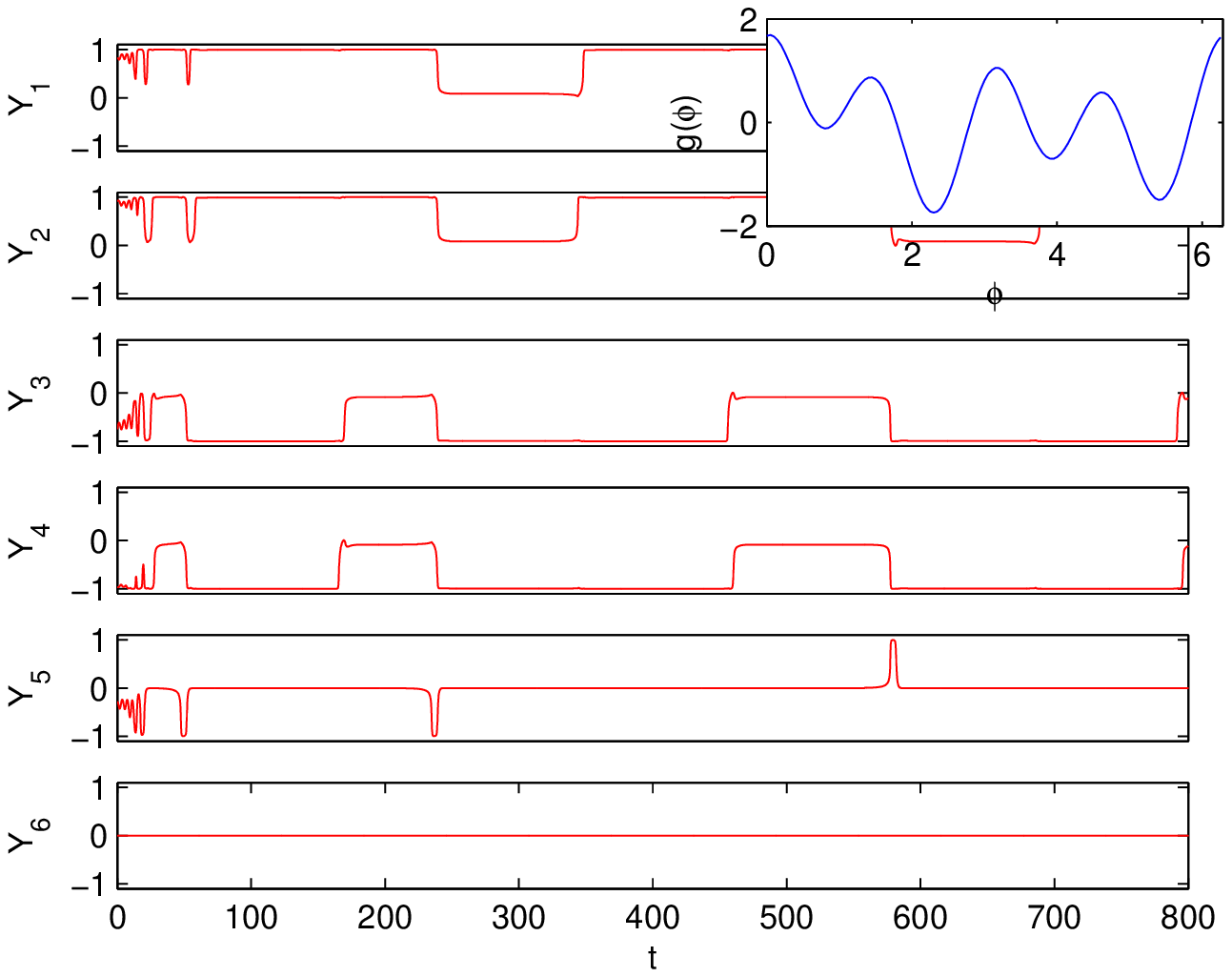}
~~\includegraphics[width=8cm]{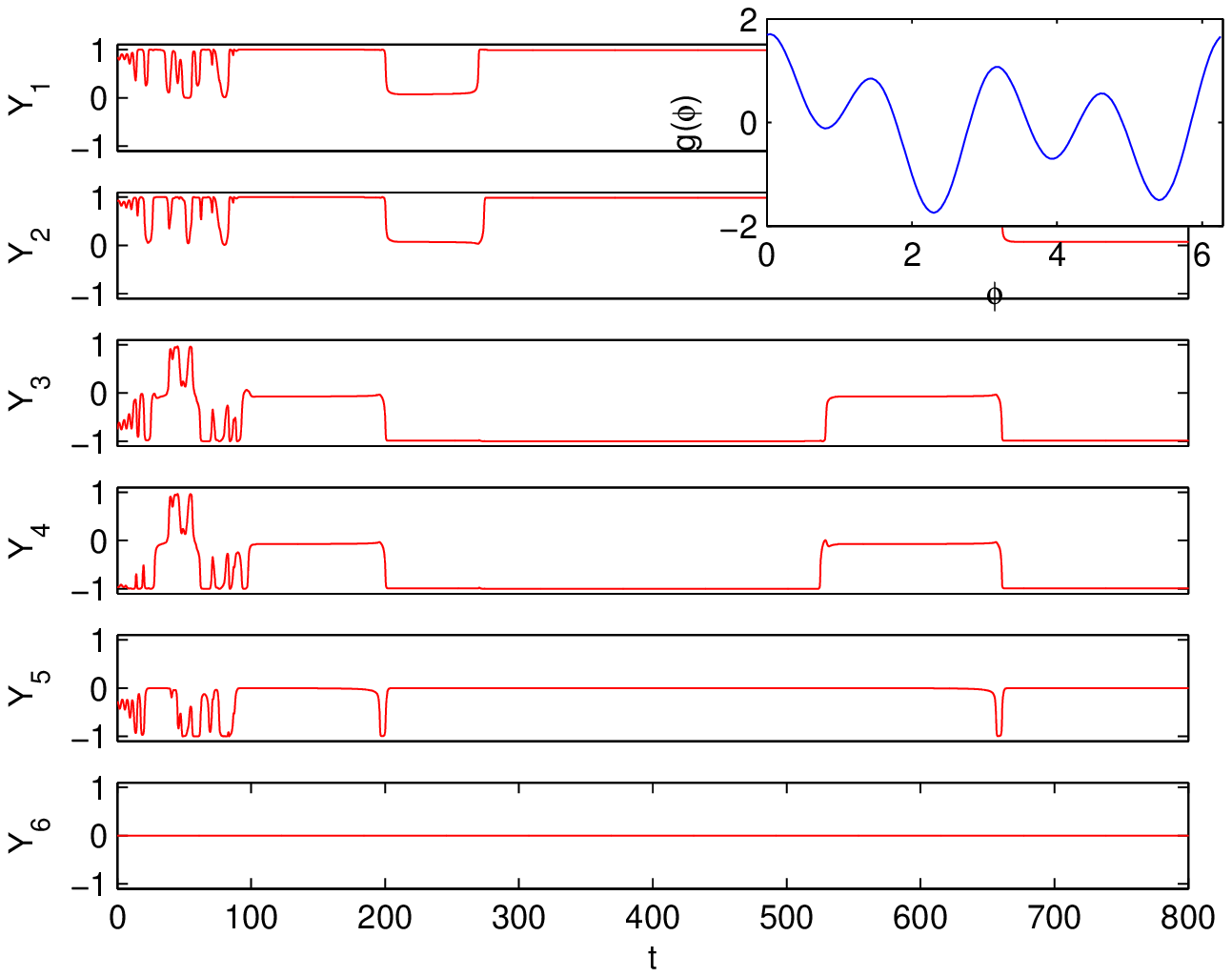}}
\caption{\label{fig:examplecase12}Left: The inset shows for Case 1 the coupling function $g(\varphi)$ and a timeseries representing the phase difference of the $k$th oscillator relative to the 6th oscillator; $Y_k= \sin(\varphi_k-\varphi_6)$ as a function of $t$. Observe that the six oscillators synchronize into three clusters for most of the time, but there are short times when the clusters break along a connection. Right: Although the cycle in Case 2 has an additional unstable direction, the trajectory still appears to approach a heteroclinic cycle between three symmetrically related states. In both cases, i.i.d. white noise of amplitude $10^{-12}$ is added to each component of the ODE.}%
\end{figure}%

The parameters listed as Cases 0, 1 and 2 in Table~\ref{tab:3clustereg} all yield three-cluster states of type $(2,2,2)$ with varying numbers of positive transverse eigenvalues. While Case 0 gives a stable cluster the dynamics for Cases 1 and 2 are more subtle. As illustrated in Figure~\ref{fig:examplecase12}, a randomly chosen initial condition evolves towards a heteroclinic cycle that connects three symmetrically related periodic cluster states within the same invariant subspace. 

Note that any phase non-degenerate $(2,2,2)$-clustered periodic orbit will have a representative that exists within the subspace $I_0$ where
$$
I_0=\{(\phi,\phi,\psi,\psi,\eta,\eta)~:~\phi,\psi,\eta\in\T\}.
$$
Note moreover that this $(S_2)^3$ invariant subspace will contain six distinct clustered periodic orbits given by cyclically permuting the phases of the clusters, due to the clusters being inequivalent.

\begin{figure}
\centerline{\includegraphics[width=0.5\columnwidth,clip=]{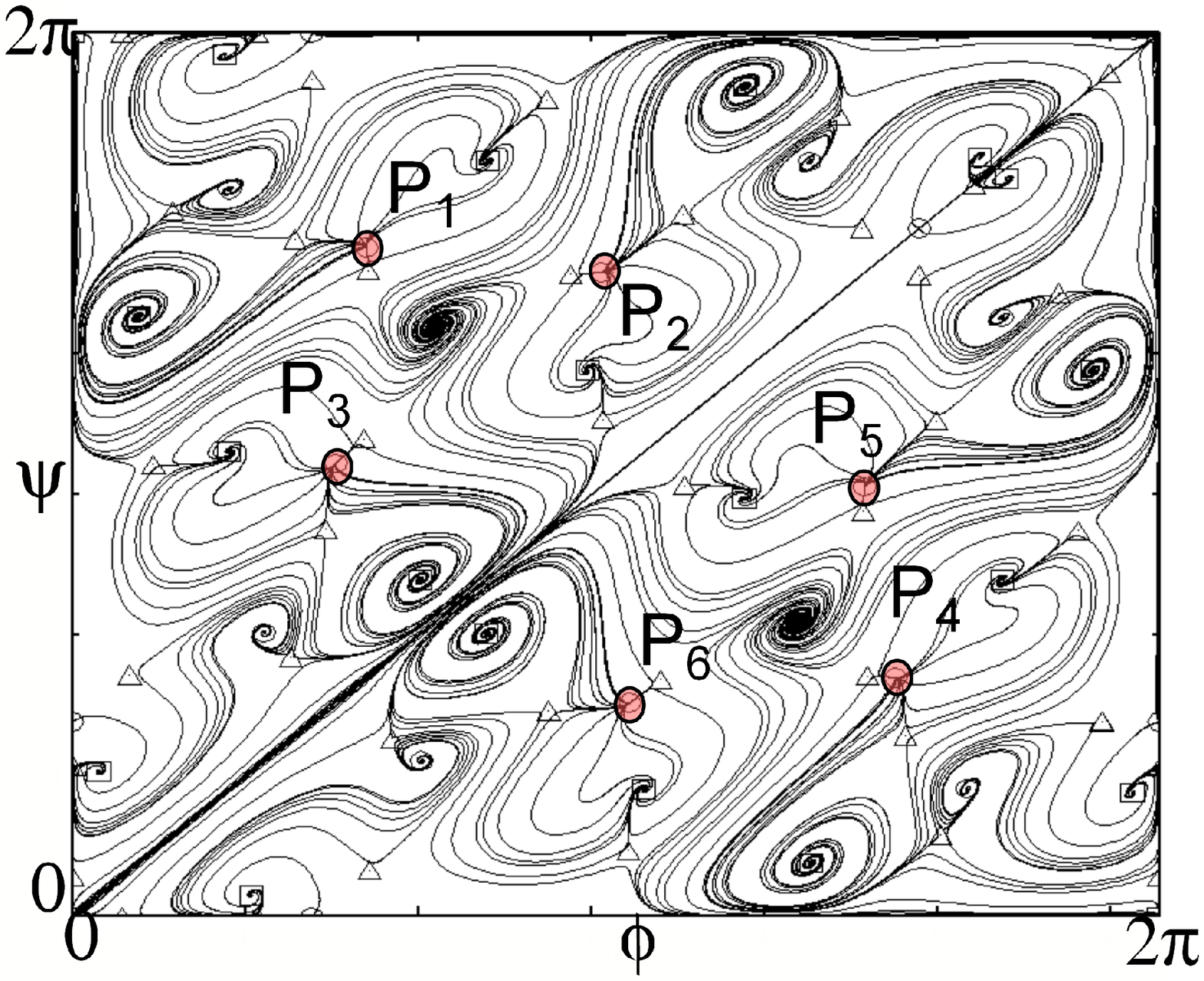}~~\includegraphics[width=0.41\columnwidth,clip=]{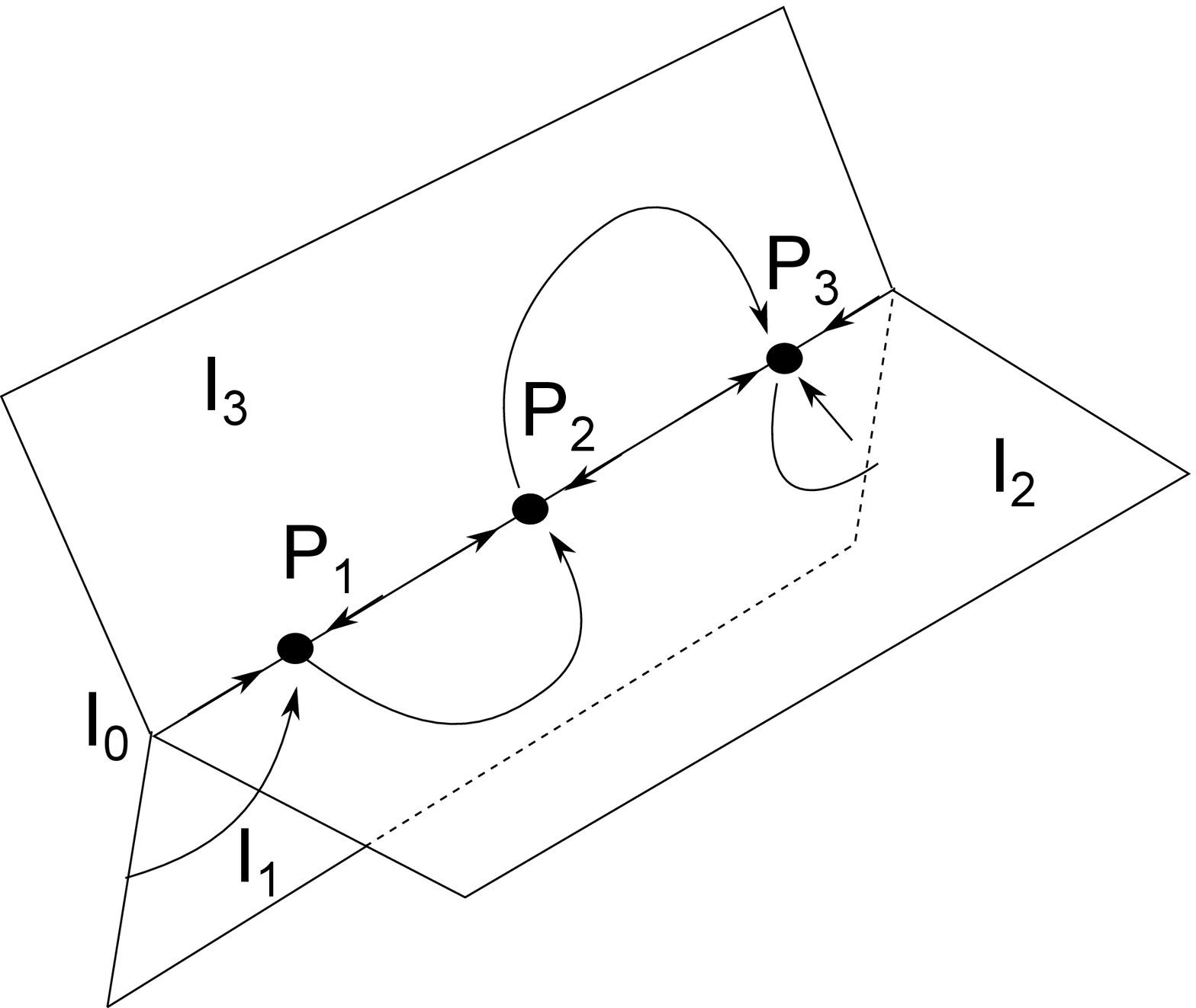}}
\caption{Left: Dynamics within the $(2,2,2)$-subspace $I_0$ showing the numerically calculated phase portrait in the $(\phi,\psi)$-plane for $\eta=0$ for Case 1; note that the projected vector field is $2\pi$ periodic in both directions. The squares, triangles and circles respectively indicate sources, saddles and sinks. The larger grey circles indicate the location of the states $P_{i}$ corresponding to inequivalent clusters within $I_0$. Right: Schematic diagram showing the heteroclinic cycle between six distinct $(2,2,2)$-cluster states for six oscillators for Case 1. The $(2,2,2)$-cluster states $P_i$ for $i=1,2,3$ within $I_0$ are connected via the fixed point subspaces $I_i$ for $i=1,2,3$ and their intersection $I_0$ while a similar arrangement connects $P_i$ for $i=4,5,6$. 
\label{fig:222clusterdiag}}%
\end{figure}

More specifically, suppose there is a point on a phase non-degenerate periodic $(2,2,2)$-cluster
$$
P_1=(0,0,\alpha,\alpha,\beta,\beta)
$$
for some $\alpha,\beta\in\T$. Without loss of generality one can choose $0<\alpha<\beta<2\pi$(in fact, one can assume that $0<2\alpha\leq \beta<2\pi-\alpha$) and phase non-degeneracy means that $\beta\neq 2\alpha$. As a consequence, $P_2$ and $P_3$ are also points on periodic $(2,2,2)$-clusters, where 
\begin{eqnarray*}
P_2&=&(0,0,\beta-\alpha,\beta-\alpha,\beta,2\pi-\alpha,2\pi-\alpha),\\
P_3&=&(0,0,2\pi-\beta,2\pi-\beta,2\pi+\alpha-\beta,2\pi+\alpha-\beta)
\end{eqnarray*}
and although these are also within $I_0$, the phase non-degeneracy means that they are distinct points; there are three more that are in the same subspace which we write as
\begin{eqnarray*}
P_4&=&(0,0,2\pi-\alpha,2\pi-\alpha,\beta-\alpha,\beta-\alpha),\\
P_5&=&(0,0,2\pi+\alpha-\beta,2\pi+\alpha-\beta,2\pi-\beta,2\pi-\beta),\\
P_6&=&(0,0,\beta,\beta,\alpha,\alpha).
\end{eqnarray*}
The relative location of these six equilibria can be seen in Figure~\ref{fig:222clusterdiag} (left) calculated using xppaut \cite{ermentrout}. For the coupling function in Case 1, Table~\ref{tab:3clustereg} reveals that each of the $P_i$ has a single positive transverse eigenvalue corresponding to instability of one of the clusters but is otherwise stable. The unstable manifold will therefore be contained within a fixed point subspace of symmetry $(S_2)^2$ where one of the clusters is broken, but the numerical results in Figure~\ref{fig:examplecase12} indicate that this unstable manifold is within the stable manifold of one of another $P_i$. If we write
\begin{itemize}
\item $I_{1}=\{(\phi,\xi,\psi,\psi,\eta,\eta)~:~\phi,\psi,\xi,\eta\in\T\}$
\item $I_{2}=\{(\phi,\phi,\psi,\psi,\eta,\xi)~:~\phi,\psi,\xi,\eta\in\T\}$
\item $I_{3}=\{(\phi,\phi,\psi,\xi,\eta,\eta)~:~\phi,\psi,\xi,\eta\in\T\}$
\end{itemize}
then one can verify that there will be a sequence of connections (heteroclinic orbit) (a) from $P_1$ to $P_2$ that is transverse within $I_{2}$, (b) from $P_2$ to $P_3$ that is transverse within $I_{3}$ and (c) from $P_3$ to $P_1$ that is transverse within $I_{1}$ as show in Figure~\ref{fig:222clusterdiag}(right)
\begin{equation}
P_1\stackrel{I_{2}}{\longrightarrow} P_2\stackrel{I_{3}}{\longrightarrow} P_3\stackrel{I_{1}}{\longrightarrow}P_1.
\label{eq:heteroclinic cycle between 3_cluster state}
\end{equation} 
Similarly there is a symmetrically related heteroclinic cycle
\begin{equation}
P_4\stackrel{I_{2}}{\longrightarrow} P_5\stackrel{I_{3}}{\longrightarrow} P_6\stackrel{I_{1}}{\longrightarrow}P_4
\label{eq:heteroclinic cycle_b}
\end{equation}
that connects the remaining equilibria within $I_0$.

Turning to Case 2 in Figure~\ref{fig:examplecase12} we note that there are two transversely unstable directions from each of the $P_i$ in Figure~\ref{fig:222clusterdiag} and ``accordingly''  a continuum of directions by which the trajectory can leave a neighbourhood of the $P_i$. These need no longer be within any of the invariant subspaces $I_i$. As can be seen for this case, the resulting dynamics ``nonetheless'' seems to return repeatedly to a cycle between the $P_i$ suggesting that the cycle is a Milnor attractor.

\section{Discussion}
\label{sec:conclusion}

The results in Section~\ref{sec:3cluster} (concerning clustering behavior and bifurcations) apply to systems with any number of coupled oscillators, regardless of the size of each cluster. However, these results are essentially local in phase space. More precisely, a given coupling function may admit a variety of different cluster states of varying stability, and there may be constraints on the possible cluster states and/or their stabilities. It would be interesting to understand the nature of  such constraints, but we leave this for future work.

By considering properties of three-cluster states with equal sized but inequivalent clusters, we find in Section~\ref{sec:3cluster} a new type of robust heteroclinic attractor for $N=6$ oscillators. Our results on the existence of robust connections for these heteroclinic attractors still rely on numerical observation of robust connections - it is a challenge to characterise coupling functions that give rise to such cycles in a more analytical (or geometric) manner. This does not seem to be an easy task, even if one restricts to cycles with one-dimensional unstable manifolds, i.e. between states that have clusterings consisting only of pairings.

Finally, we briefly state a generalization of Corollary~\ref{cor:couplingfunction} to multi-cluster states of arbitrary size.

\begin{theorem}
\label{thm:couplingfunctiongen}
Suppose that $g$ is such that there is a periodic $M$-cluster state $(\psi_1,\ldots,\psi_M)$ with non-trivial clusters of size $m_1,\ldots,m_M$, such that the clusters are tangentially stable. Then there is a parametrized family of coupling functions
$$
g_{r}(\varphi)=g(\varphi)+r h(\varphi)
$$
and with real parameter $r\in \R$ and parameter values $r_1 \leq \ldots \leq r_M$ such that for all values of the parameter $r$ a nearby cluster state exists, is still tangentially stable and moreover:
\begin{enumerate}
\item If $r<r_1$ then all clusters are unstable.
\item If $r_k<r<r_{k+1}$ then precisely $k$ of the clusters are stable.
\item If $r_M<r$ then all clusters are stable.
\end{enumerate}
\end{theorem}

\proof
The proof is similar to that of Corollary~\ref{cor:couplingfunction}; we make use of the fact that the transverse exponent of the $k$th cluster can be written in the form
$$
\lambda_k^{tran}=\frac{1}{N} \left[m_kg'_0+\sum_{l=1,l\neq k}^{M} m_lg'_{kl}\right].
$$
The nontrivial assumption of the clusters mean one can choose a compactly supported perturbation $h$ with $h(0)=0$ and $h'(0)=-1$ so that $g_{r,0}$, $g_{r,kl}$ and $g'_{r,kl}$ are independent of $r$ while $g'_{r,0}=g'_0-r$.
\qed

~

\subsection*{Acknowledgements}

We thank Christian Bick, Hiroshi Kori and Mike Field for some very interesting conversations with regard to this work.

\bibliographystyle{plain}


\end{document}